\documentclass[krantz1,ChapterTOCs]{krantz-grg} 
\usepackage{fixltx2e,fix-cm}
\usepackage{amssymb}
\usepackage{amsmath}
\usepackage{amsfonts}
\usepackage{amssymb}
\usepackage{amsthm}
\usepackage{subfigure}
\usepackage{graphicx}
\usepackage{color}
\usepackage{enumerate}
\usepackage{lineno}
\usepackage{graphicx}
\usepackage{makeidx}
\usepackage{multicol}
\usepackage{lipsum}
\usepackage{pinlabel}
\usepackage{verbatim}

\newtheorem{theorem}{Theorem}[section]
\newtheorem{proposition}[theorem]{Proposition}

\theoremstyle{definition}
\newtheorem{definition}[theorem]{Definition}

\newtheorem{conjecture}[theorem]{Conjecture}

\theoremstyle{remark}

\numberwithin{equation}{section}

\AtBeginDocument{\renewcommand{\bibname}{Resource List}}
\HeadingsChapterSection

\makeindex  

\usepackage{bm,url}
\newcommand\fpq{\phi_{p,q}}
\newcommand\om{\omega}
\newcommand\Om{\Omega}
\newcommand\qq{\qquad}
\newcommand\q{\quad}
\newcommand\rc{random-cluster}
\newcommand\RR{{\mathbb R}}
\newcommand\s{\sigma}
\newcommand\ZI{Z_{\text{\rm I}}}
\newcommand\HI{H_{\text{\rm I}}}
\newcommand\HP{H_{\text{\rm P}}}
\newcommand\ZP{Z_{\text{\rm P}}}
\newcommand\ZRC{Z_{\text{\rm RC}}}
\newcommand\la{\lambda}
\newcommand\ra{\rangle}
\newcommand\be{\beta}
\newcommand\de{\delta}
\newcommand\lra{\leftrightarrow}

\newcommand\PP{{\mathbb P}}
\newcommand\EE{{\mathbb E}}
\newcommand\ZZ{{\mathbb Z}}
\newcommand\LL{{\mathbb L}}

\newcommand\Si{\Sigma}
\newcommand\oo{\infty}
\newcommand\lest{\le_{\text{\rm st}}}
\newcommand\gest{\ge_{\text{\rm st}}}
\newcommand\Omf{\Om_{\text{\rm for}}}
\newcommand\Omst{\Om_{\text{\rm st}}}
\newcommand\Omcs{\Om_{\text{\rm cs}}}
\newcommand\UST{\text{\rm UST}}
\newcommand\UCS{\text{\rm UCS}}
\newcommand\USF{\text{\rm USF}}
\newcommand\sm{\setminus}
\newcommand\comp[1]{\overline{#1}}
\newcommand\ssquare{\,\square\,}
\newcommand\ov{\vec}
\newcommand\es{\varnothing}
\newcommand\boldm{{\mathbf m}}
\newcommand\boldj{{\mathbf J}}
\newcommand\boldv{{\mathbf v}}
\newcommand\boldh{{\mathbf h}}
\newcommand\bzero{{\mathbf 0}}
\newcommand\bl{\beta} 
\newcommand\bP{{\text{\bf P}}}
\newcommand\pim{P}
\newcommand\resp{respectively}
\newcommand\lc{\lambda_{\text{\rm c}}}
\newcommand\tmax{\theta_{\text{\rm max}}}
\renewcommand\Pr{\PP}

\newcounter{mycount}
\newenvironment{romlist}{\begin{list}{\rm(\roman{mycount})}%
   {\usecounter{mycount}\labelwidth=1cm\itemsep 0pt}}{\end{list}}
\newenvironment{numlist}{\begin{list}{\arabic{mycount}.}%
   {\usecounter{mycount}\labelwidth=1cm\itemsep 0pt}}{\end{list}}
\newenvironment{letlist}{\begin{list}{\rm(\alph{mycount})}%
   {\usecounter{mycount}\leftmargin=1cm\labelwidth=1cm\itemsep 0pt}}{\end{list}}

\begin{document}
\frontmatter

\mainmatter

\chapter*{}

\chapterauthor{Geoffrey Grimmett}{University of Cambridge\\
Centre for Mathematical Sciences,\\ Wilberforce Road,\\ Cambridge CB3 0WB, UK\\
\texttt{g.r.grimmett@statslab.cam.ac.uk}\\}

\chapter{The Potts and random-cluster models}

\section{Synopsis}

\begin{itemize}
\item Ising,  Potts, and random-cluster models
\item Basic properties of random-cluster measures
\item Flow polynomial
\item Limit as $q \downarrow 0$
\item Zero-temperature limit
\item Asymptotics of the Tutte polynomial on the complete graph

\end{itemize}

\section{Introduction}

The four principal elements in this chapter are the Ising model of 1925, the Tutte polynomial
of 1947, the Potts model of 1952, and the random-cluster model of 1972. 

The \emph{Ising 
model}\footnote{The Ising model was proposed to Ising by Lenz, and the Potts model to Potts
by Domb.}
\cite{Ising} is the fundamental model for the ferromagnet, and it has generated enormous
interest and activity in mathematics and physics over the intervening decades. The \emph{Potts
model} \cite{MR0047571} extends the number of local states of the Ising model from $2$ 
to a general number $q$. The \emph{random-cluster model} of Fortuin and Kasteleyn \cite {MR0359655}
provides an overarching framework for the Ising/Potts models that
incorporates percolation and electrical networks, together with certain other processes.
The common aspect of importance for these three systems is the singularity
that occurs at points of phase transition.

Whereas these three processes originated in mathematical physics, the 
\emph{Tutte polynomial} \cite{MR0018406} is an object from combinatorics, and it
encapsulates a number of significant features of a graph or matroid. It turns
out that the Tutte polynomial is  equal (subject to a change of variables) to
the partition function of the random-cluster model, and therefore to that of
the Potts model. This connection is not a coincidence since both the Tutte and
random-cluster functions arose in independent analyses of local graph operations
such as deletion and contraction.

The Tutte polynomial originated in Tutte's exploration of 
deletion and contraction on a finite graph. The random-cluster
model originated similarly in Kasteleyn's observation that the Ising,
Potts, and percolation models, and also electrical networks, have
a property of invariance under series and parallel operations on edges.

Combinatorial theory and statistical mechanics are areas of science 
which have much in common,
while retaining their distinctive characteristics. Statistical mechanics
is mostly concerned with the structure of phases and of singularities,
and has developed appropriate methodology and language. 
Although, in principle, the properties of a physical model
are encoded entirely within its partition function, the extraction of such properties
is often challenging and hinges frequently on other factors such as the
nature of the underlying graph.

The connection between the Tutte polynomial and statistical
mechanics is summarised in this chapter, as follows. We aim: (i)
to give a clear formulation of the relevant models, (ii) to explain
the connection between their partition functions and the Tutte polynomial,
(iii) to present some of the basic properties of the \rc\  model
that are directly contingent on the partition function, and (iv) to present a selection
of open problems concerning Potts and
\rc\ models that may be related to the Tutte polynomial.

\section {Probabilistic Models from Physics}

\subsection{The Ising and Potts ferromagnets}

The Ising model for ferromagnetism was analysed in one dimension in 
Ising's thesis and 1925 paper \cite{Ising}.
It modelled the following physical experiment.
A piece of iron is placed in a magnetic
field, with an intensity that is increased from zero to a maximum, and then
reduced to zero.  The iron
retains some residual magnetisation if and only if
the temperature is sufficiently low, and the
critical temperature for this phenomenon is called the \emph{Curie point}. 

The Ising model may be summarised as follows.
Suppose that particles are placed at the vertices of
a graph embedded  in a Euclidean space. Each particle may be in either
of two states: spin `up' or spin `down'. Spin-values are chosen at random according to
a certain probability measure governed by interactions between neighbouring particles.
This measure is described as follows.

Let $G=(V,E)$ be a finite, simple graph.
Each vertex $v\in V$ is occupied by a particle with a random spin.
Since spins are assumed to come in two basic types, we take as sample space
the set $\Sigma =\{ -1,+1\}^V$ of vectors $\s=(\s_x: x\in V)$ with entries $\pm 1$. 

\begin{definition} \label{def:isi}
Let  $\be\in(0,\oo)$ and $h\in\RR$. The Ising probability measure $\la_{\be,h}$ 
on $G$ is given by
$$
\la_{\be,h} (\s )=\frac{1}{\ZI}e^{-\be \HI(\s )},
\qq\s\in\Sigma ,
$$
where the \emph{partition function} $\ZI$ and the \emph{Hamiltonian} $\HI$ 
are given by
$$
\ZI(\be,h)= \sum_{\s\in\Sigma} e^{ -\be \HI(\s )},\qq
\HI(\s )=-\sum_{e=\langle x,y\rangle\in E} \s_x\s_y-h\sum_{x\in V}\s_x.
$$
\end{definition}

The parameter $\be$ represents the reciprocal $1/T$ of
temperature, and $h$ is the \emph{external field}.
For reasons of simplicity, we assume usually that $h=0$, and
we write $\la_{\be} = \la_{\be,0}$. 
It is usual to include also an edge-interaction
$J$, which we have chosen to wrap within the parameter $\be$.

The Ising model has two admissible local spin-values, and a very rich theory.
In his 1952 paper \cite{MR0047571}, Potts developed an extension of the Ising model
to a general number of spin-values.

Let $q$ be an integer satisfying  $q\ge 2$, and consider the sample space
$\Sigma =\{ 1,2,\ldots ,q\}^V$. Each vertex of $G$ may now be
in any of $q$ states. 

\begin{definition}\label{def:pot}
Let $\be\in(0,\oo)$ and $q\in \{2,3,\dots\}$. 
The Potts probability measure on $G$ is 
given by
$$
\pi_{\be,q}(\s )=\frac{1}{\ZP}e^{ -\be  \HP(\s )},
\qq\s\in\Sigma,
$$
where the Potts partition function $\ZP$ and Hamiltonian $\HP$ are given by
$$
\ZP(\be,q)= \sum_{\s\in\Sigma} e^{ -\be \HP(\s )},\qq
\HP(\s )=-\sum_{e=\langle x,y\ra} \de_{\s_x,\s_y},
$$
and $\de_{u,v}$ is the Kronecker delta.
\end{definition}

When $q=2$, we have that 
$$
\de_{\s_x,\s_y}=\tfrac12(1+\s_x\s_y),
$$
from which it follows that the $q=2$ Potts model is simply the Ising
model with $\be$ replaced by $\frac12\be$.

In a more general definition, one may include a non-zero external field $h$ and
a vector $\boldj$ of edge-parameters. See Section \ref{sec: more}.

\subsection{The random-cluster model}

The random-cluster model was formulated in a series of papers 
\cite{ MR0378660,MR0432137,MR0359655} by Fortuin
and Kasteleyn. It is described next, and its relationship to the Potts model
is explained in Section \ref{sec:coupling}.

Let $G=(V,E)$ be a finite, simple graph. 
The relevant state space is the set
$\Omega=\{0,1\}^E$ of vectors $\omega=(\omega(e):e\in E)$
with $0/1$ entries.
An edge $e$ as said to be \emph{open} in $\omega\in \Om$ if
$\omega(e)=1$, and \emph{closed} if $\omega(e)=0$.
For $\omega\in\Omega$, let $\eta(\omega)=\{e\in E: \omega(e)=1\}$
denote the set of open edges, and let $k(\omega)$ be the number
of connected components (or `open clusters') of the graph
$(V,\eta(\omega))$; the count $k(\omega)$ includes the number of isolated
vertices.

\begin{definition}\label{def:rcm}
Let $p\in(0,1)$ and $q\in(0,\oo)$.
The \emph{random-cluster measure}  $\fpq$ on $G$ is given by
$$
\fpq(\omega) = \frac 1\ZRC \biggl\{\prod_{e\in E} p^{\om (e)} (1-p)^{1-\om
(e)}\biggr\} q^{k(\om )},\qq\om\in\Om,
$$
where the partition function $\ZRC$ is given by
\begin{equation}\label{eq:rcpart}
\ZRC(p,q)=\sum_{\om\in\Om}\biggl\{\prod_{e\in E} p^{\om (e)} (1-p)^{1-\om
(e)}\biggr\} q^{k(\om )}.
\end{equation}
\end{definition}

The most important values of $q$ are arguably the positive
integers. When $q=1$, $\phi_p := \phi_{p,1}$ is a product measure, and the
words `percolation' and `random graph' are often used
in this context, see \cite{MR1707339,
 MR1782847}.
 The \rc\ model with
$q\in\{2,3,\dots\}$ corresponds, as sketched in
the next section, to the Potts model with $q$ local states.

See \cite{MR2243761} for the general theory of the random-cluster model.

\subsection{Coupling of the Potts and random-cluster measures}\label{sec:coupling}

Let $q\in\{2,3,\dots\}$, $p\in (0,1)$,
and let $G=(V,E)$ be a finite, simple graph.  
We consider the product sample space $\Sigma\times\Om$ where
$\Sigma=\{1,2,\dots,q\}^V$ and $\Om=\{0,1\}^E$ as above.
Let $\mu$ be the probability measure on $\Sigma\times\Om$ given by
$$
\mu(\s,\om) \propto \psi(\s)\phi_p(\om) 1_F(\s,\om),\qq
(\s,\om)\in\Si\times\Om,
$$
where $\psi$ is uniform measure on $\Si$,
$\phi_p=\phi_{p,1}$ is product measure
on $\Om$ with density $p$, and $1_F$ is the indicator function
of the event that $\s$ is constant on each open cluster of $\om$, that is,
$$
F=\bigl\{ (\s,\om)\in\Si\times\Om : \text{$\s_x=\s_y$ for every
$e=\langle x,y\ra$ satisfying $\om(e)=1$}\bigr\}.
$$
The measure $\mu$ may be viewed as the product measure
$\psi\times\phi_p$ conditioned on the event $F$.

\begin{theorem}
Let $q\in\{2,3,\dots\}$ and $p\in (0,1)$.
\begin{letlist}
\item 
\emph{Marginal on $\Sigma$}. The first marginal of $\mu$ (on $\Si$)
is the Potts measure $\pi_{\be,q}$
where $p=1-e^{-\be}$. 

\item
\emph{Marginal on $\Om$}. The second marginal of $\mu$ (on $\Om$) is
the \rc\ measure $\fpq$.

\item
\emph{Conditional measures}. 
\begin{romlist}
\item
Given $\om\in\Om$, the
conditional measure on $\Si$ is obtained by putting (uniformly)
random spins on entire clusters of $\om$. 
These spins are constant on given clusters, and are independent between
clusters. 
\item
Given $\s\in\Si$, the conditional measure on $\Om$ is obtained as follows.
For $e=\langle x,y\ra\in E$, we set $\om (e)=0$ if $\s_x\ne\s_y$, and otherwise $\om (e)=1$ with
probability $p$ (independently of other edges).
\end{romlist}

\item \emph{Partition functions.} We have that
\begin{equation}\label{eq:partition}
\ZRC(p,q)=e^{-\be |E|}\ZP(\be,q).
\end{equation}
\end{letlist}
\end{theorem}

This coupling may be used  
to show that correlations
in Potts models correspond to connection probabilities in random-cluster models.
In this way, one may harness methods of stochastic geometry in order to
understand the correlation structure of the Potts system.
The basic theorem of this type is Theorem \ref{corrconn}, following.

The `two-point correlation function' of the
Potts measure $\pi_{\be,q}$ on the finite graph 
$G=(V,E)$ is the function
\begin{equation}
\label{eq:tau}
\tau_{\be,q}(x,y):=\pi_{\be,q}(\s_x=\s_y)-\frac{1}{q},\qq x,y\in V.
\end{equation}
Let $\{x\lra y\}$ be the event of $\Omega$ on which
there exists an open path joining vertex $x$ to vertex $y$. 
The `two-point connectivity function' of the random-cluster measure
$\phi_{p,q}$ is the function $\fpq (x\lra y)$ for $x,y\in V$, 
that is, the probability that $x$ and $y$
are joined by a path of open edges.
It turns out that these two-point functions
the same up to a constant factor.

\begin{theorem} [Correlation/connection theorem]  \label{corrconn}
If $q\in\{ 2,3,\ldots\}$ and $p=1-e^{-\be}\in (0,1)$, then
\begin{equation}\label{eq:corrconn}
\tau_{\be,q}(x,y)=(1-q^{-1})\fpq (x\lra y),\qq x,y\in V.
\end{equation}
\end{theorem}

The Potts models considered above have zero external field. Some 
complications arise when an external field is added; see Section
\ref{sec: more}.

\subsection{Partition functions and the Tutte polynomial}

The Potts and \rc\ partition functions may be viewed as
evaluations of rank-generating functions and the Tutte polynomial, defined as follows.
The (Whitney) \emph{rank-generating function} of $G=(V,E)$ is the function
$$
W_G(u,v)=\sum_{A\subseteq E} u^{r(A)} v^{c(A)},\qq u,v\in\RR,
$$
where 
$$
r(A)=|V|-k(A), \qq c(A)=|A|-|V|+k(A),
$$
are the rank and co-rank of $G'=(V,A)$, \resp.
Here, $k(A)$
denotes the number of components of $G'$.
The Tutte polynomial may be written as
$$
T_G(u,v)= (u-1)^{|V|-1} W_G\bigl( (u-1)^{-1}, v-1\bigr).
$$

\begin{theorem}\label{thm:tutte-rcm}
Let $p\in (0,1)$, $q\in (0,\oo)$, and
$$
u-1=\frac{q(1-p)}p,\qq v-1=\frac p{1-p}.
$$
\begin{letlist}
\item 
The partition function $\ZRC=\ZRC(G)$ of the random-cluster
measure on $G$ with parameters $p$, $q$ satisfies
\begin{align*}
\ZRC(G)=\left[(u-1)(v-1)^{|V|} v^{-|E|}\right] T_G(u-1,v-1).
\end{align*}

\item If $q \in\{2,3,\dots\}$ and $p=1-e^{-\be}$,
the partition function of the $q$-state Potts model satisfies
$$
\ZP(\be,q) =\left[(u-1)(v-1)^{|V|}\right] T_G(u-1,v-1).
$$
\end{letlist}
\end{theorem}

\subsection{Potts extensions}\label{sec: more}

The Potts model of Definition \ref{def:pot} is in its simplest form in three senses:
(i) each edge plays an equal (deterministic) role, (ii) the external field  satisfies $h=0$,
and (iii) the model is \emph{ferromagnetic}.  More generally, one may
consider the partition function
$\ZP(\be,\boldj,h) = \sum_{\s\in\Si} e^{-\be \HP(\s)}$ where the 
Hamiltonian is given by
\begin{equation}\label{eq:genpot}
\HP(\s )=-\sum_{e=\langle x,y\ra} J_e\de_{\s_x,\s_y} - \sum_{j=1}^q
\sum_{x \in V} h_j\de_{\s_x,j}.
\end{equation}
Here, $\boldj=(J_e: e \in E)$ is a family of edge-parameters assumed
to satisfy $J_e\ne 0$,
and  $\boldh=(h_x: x\in V)$
is a vector of external
fields, one for each possible local spin $1,2,\dots,q$. The model is termed \emph{ferromagnetic}
if $J_e > 0$ for all $e \in E$, and \emph{purely antiferromagnetic} if $J_e < 0$ for all $e\in E$.
The  general Potts partition function of \eqref{eq:genpot} poses some new difficulties.

Assume first that $\boldh=\bzero$.
The associated \rc\ formula yields a function $\phi_{\bold p,q}$ where
$p_e=1-e^{-\be J_e}$. If $J_e<0$ for some $e$, this does not define a probability measure.
In addition, the Potts model does not satisfy the range of correlation inequalities
that hold in the ferromagnetic case. 
On the other hand, Theorem \ref{thm:tutte-rcm}(b) is easily extended 
for general $\boldj$ to a 
multivariate Tutte polynomial on  $G=(V,E)$ which
may be written in the form
\begin{equation}\label{eq:mtutte}
T_G(q,\boldv) = \sum_{A \subseteq E} q^{k(A)}\prod_{e\in A}v_e,
\end{equation}
where $q$ and $\boldv = (v_e: e \in E)$ are viewed as parameters.
See \cite{MR2832375,MR2187739} for recent 
accounts.

The $J_e$ may themselves be random, in which case the model is
termed \emph{quenched}, in contrast to the \emph{annealed} case
in which one averages initially over the random  $J_e$. If the
probability distribution
of the $J_e$ allocates mass to both positive and 
negative values, the system is a `spin glass'. See \cite{NewS}.

When $\boldh \ne \bzero$, a form of the Tutte--Potts correspondence 
may be found in \cite{MR1757955},
where positive association and infinite-volume limits are explored, 
and also  in a slightly more general setting
in \cite{MR2832375} (see also \cite{MR1703438}).

\section{Basic Properties of Random-Cluster Measures}

This section includes some of the basic properties of a \rc\ measure
on the finite graph $G=(V,E)$.

\subsection{Stochastic ordering}

The state space $\Om=\{0,1\}^E$ is a partially ordered set with partial order
given by: $\om_1\le\om_2$ if $\om_1(e)\le\om_2(e)$ for all $e \in E$. 
This partial order induces partial orderings on
the spaces of associated functions and measures,
and this is extremely useful in the
analysis of Potts and \rc\ models.

\begin{definition}\label{def:mon}
\mbox{\hfil}
\begin{letlist}
\item
A random
variable $f:\Om\to\RR$ is called \emph{increasing} if
$f(\om_1)\le f(\om_2)$ whenever $\om_1\le \om_2$.

\item
An event $A\subseteq \Om$ is called \emph{increasing} if its indicator
function $1_A$ is increasing. 

\item
Given two probability measures
$\mu_1$, $\mu_2$ on $\Om$, we write $\mu_1\lest \mu_2$, and
say that $\mu_1$ is stochastically smaller than $\mu_2$,
if $\mu_1(f)\le \mu_2(f)$ for all 
increasing random variables $f$ on $\Om$.
\end{letlist}
\end{definition}

Arguably the most useful approach to stochastic ordering is
due to Holley.  We obtain the following comparison inequalities
as corollaries of Holley's inequality, see \cite{MR0341552} and \cite[Thm 2.1]{MR2243761}.

\begin{theorem}[Comparison inequalities]\label{dom}
It is the case that
\begin{align*}
\phi_{p',q'}&\lest\phi_{p,q}\qq\text{if}\q q'\geq q,\ q'\geq 1,
\text{ and } p'\leq p,\\
\phi_{p',q'}&\gest\phi_{p,q}\qq\text{if}\q q'\geq q,\ q'\geq 1,
\text{ and } \frac{p'}{q'(1-p')} \geq \frac p{q(1-p)}.
\end{align*}
\end{theorem}

\subsection{Positive association}

Holley's inequality admits a neat proof of the FKG inequality of \cite{MR0309498}.
This amounts to the following in the case of \rc\ measures.

\begin{theorem}[Positive association] \label{thm:pa}
Let $p\in(0,1)$ and $q\in[1,\oo)$. If $f$
and $g$ are increasing functions on $\Om$, then
$$\fpq (fg)\geq \fpq(f)\fpq (g).
$$
\end{theorem}

Specializing to indicator functions $f=1_A$, $g=1_B$, we obtain that
$$
\phi_{p,q} (A\cap B)\geq\phi_{p,q} (A)\phi_{p,q} (B)\qq\hbox{for
increasing events $A,B$},
$$
whenever $q\geq 1$. Positive association is generally false when $0<q<1$.

\section{The Limit as $q\downarrow 0$}

\subsection{UST, USF, and USC}

Some interesting limits with combinatorial flavours 
arise from consideration of $\fpq$ as $q\downarrow 0$. 
Write $\Omf$, $\Omst$, $\Omcs$ for the subsets of $\Om$ containing all
forests, spanning trees, and connected subgraphs of the
underlying graph $G$, respectively, and
write $\USF$, $\UST$,  and $\UCS$ for the 
uniform probability measures
on the respective sets $\Omf$, $\Omst$, $\Omcs$.
An account of the following limits and their history
may be found at \cite[Thm 1.2]{MR2243761}.

\begin{theorem}
We have in the limit as $q\downarrow 0$ that{\rm:}
$$
\phi_{p,q}\Rightarrow \begin{cases} \UCS &\text{if } p=\tfrac12,\\
 \UST &\text{if } p\to 0 \text{ and } q/p\to 0,\\
 \USF &\text{if } p=q.
\end{cases}
$$
\end{theorem}

The spanning tree limit $\UST$ is especially interesting for historical and
mathematical reasons.
The \rc\ model originated in a systematic
study by Fortuin and Kasteleyn of systems of 
a certain type which satisfy certain parallel and series laws. Electrical 
networks are the best known such systems:
two resistors of resistances $r_1$ and
$r_2$ in parallel (respectively, in series) 
may be replaced by a single resistor 
with resistance
$(r_1^{-1}+r_2^{-1})^{-1}$
(respectively, $r_1+r_2$). Fortuin and Kasteleyn realized
that the electrical-network theory of a graph $G$ is
related to the limit as $q\downarrow 0$ of the \rc\ model
on $G$, where $p$ is given  by $p=\sqrt q/(1+\sqrt q)$. 
It has been known since Kirchhoff's
theorem \cite{Kirchhoff} that the electrical currents which flow
in a network may be expressed in terms of counts of spanning trees. 

The theory of the uniform spanning tree measure UST is 
beautiful in its own right, and is
linked in an important way to the emerging field of stochastic
growth processes of `stochastic L\"owner evolution' (SLE) type.
See \cite[Chap. 2]{MR2723356}.

\subsection{Negative association}

Let $E$ be
a  finite set, and let $\mu$ be a probability measure on the space $\Om=\{0,1\}^E$.
There are several concepts of negative association, of 
which we present three here.  

For $\om \in\Om$ and $F\subseteq E$, the
cylinder event $\Om_{F,\om}$ generated by $\om$ on $F$ is given by
$$
\Om_{F,\om}=\{\om'\in\Om :\om'(e)=\om(e) \text{ for } e\in F\}.
$$
For $F\subseteq E$ and an event
$A\subseteq \Om$, we say that 
$A$ \emph{is defined on} $F$ if, for all $\om\in\Om$,
we have that $\om\in A$ if and only if $\Om_{F,\om} \subseteq A$.
Let $A$ and $B$ be events in $\Om$.
We define $A \ssquare B$ to be the set of all vectors 
$\om\in\Om$ for which there
exists $F\subseteq E$ such that $\Om_{F,\om} \subseteq A$ and
$\Om_{\comp{F},\om}\subseteq B$, where 
$\comp{F}=E\sm F$. Note that the choice
of $F$ is allowed to depend on the vector $\om$. The operator
$\square$ originated in the work of van den Berg and Kesten
\cite{MR799280} on the well known BK inequality.

\begin{definition}
\mbox{\hfil}
\begin{letlist}
\item
The measure $\mu$ is 
\emph{edge negatively associated} if
$$
\mu(J_e\cap J_f) \le \mu(J_e)\mu(J_f), \qq e,f\in E,\ e\ne f,
$$
where $J_e$ is the event that $e$ is open.

\item
We call $\mu$ \emph{negatively associated} if 
$$
\mu(A \cap B) \le \mu(A)\mu(B)
$$ 
for all pairs $(A,B)$ of increasing events with the property that
there exists $F\subseteq E$ such that $A$ is defined on $F$, and $B$ is defined on 
its complement $E\sm F$.  

\item

We say that $\mu$ has the \emph{disjoint occurrence property} if
$$
\mu(A\ssquare B) \le \mu(A)\mu(B),\qq  A,B \subseteq \Om.
$$

\end{letlist}
\end{definition}
 
\begin{proposition}
We have that
\begin{align*}
&\text{$\mu$ has the disjoint occurrence property}\\
&\hskip1cm\Rightarrow\text{$\mu$ is negatively associated}\\
&\hskip1cm\Rightarrow\text{$\mu$ is edge negatively associated}.
\end{align*}
\end{proposition}

It was proved by Reimer \cite{MR1751301} that the product measures
$\phi_{p,1}$  have the
disjoint occurrence property.  The \rc\ measure $\phi_{p,q}$
cannot (generally) be edge negatively associated when $q>1$. 
It may 
be conjectured that $\phi_{p,q}$ satisfies some form of 
negative association 
when $q < 1$. Such a property would be very useful in studying \rc\ measures
when $q<1$.

In the absence of a satisfactory approach to the general case of
\rc\ measures with $q<1$, we turn next
to the issue of negative association of $\fpq$ in the limit
as $q\downarrow 0$.

\begin{conjecture}\label{2.forcs}
For any finite graph $G=(V,E)$, 
the uniform spanning forest measure $\USF$ and the 
uniform connected subgraph
measure $\UCS$ are edge negatively associated.
\end{conjecture}

A stronger version of this conjecture is 
that USF and UCS are negatively associated
in one or both of the further senses described above.
Numerical evidence for the conjecture is found in \cite{MR2060630}.

The $\UST$ measure is, in contrast, much better understood,
owing to the theory of electrical networks and, more particularly, Kirchhoff's
matrix--tree theorem \cite{Kirchhoff} and its ramifications. The following
was proved by Feder and Mihail \cite{FedM}.

\begin{theorem}
The uniform spanning tree measure $\UST$ is negatively associated.
\end{theorem}

The material in this section may be found in expanded form in \cite{MR2243761,MR2723356}.

\section{Flow Polynomial}

\subsection{Definition of the flow polynomial}

We turn $G$ into an oriented graph by allocating
an orientation to each edge $e\in E$, 
and the resulting
digraph is denoted $\ov G=(V, \ov E)$. If the edge $e=\langle u,v\ra\in E$ is 
oriented from $u$ to $v$, 
we say that $e$ \emph{leaves} $u$ and 
\emph{arrives} at $v$.

Let $q\in\{2,3,\dots\}$.
A function $f:E\to\{0,1,2,\dots,q-1\}$ is called
a \emph{mod-$q$ flow} on $\ov G$ if
\begin{equation}
\sum_{\text{\scriptsize$\ov e\in\ov E:$}\atop \text{\scriptsize$\ov e$ leaves $v$}} f(\ov e) -
\sum_{\text{\scriptsize$\ov e\in\ov E:$}\atop \text{\scriptsize$\ov e$ arrives at $v$}} f(\ov e) =
0\q\text{modulo $q$,\qq $v\in V$},
\label{8.flowdef}
\end{equation}
so that flow is conserved at every vertex.
A mod-$q$ flow $f$ is called \emph{non-zero}
if $f(\ov e)\ne 0$ for all
$\ov e\in \ov E$. 
Let $C_G(q)$ be the number of non-zero mod-$q$ flows on $\ov G$.

It is standard that $C_G(q)$ does not
depend on the orientations of the edges of $G$. 
The function $C_G(q)$, viewed as a function of $q$, is 
called the \emph{flow polynomial} of $G$. We adopt the convention
that $C_G(q)=1$ if $E=\es$.

The flow polynomial may be obtained
as evaluations of the Tutte and rank generating polynomial with
two particular parameter values, as follows,
\begin{align}
C_G(q)&=(-1)^{|E|} W_G(-1,-q)
\nonumber\\
&= (-1)^{|E|-|V|+1}T_G(0,1-q),
\qq q\in\{2,3,\dots\}, \label{8.flowpodef} 
\end{align}
and thus $C_G$ is indeed a polynomial. 
We shall later write $C(G;q)$ for
$C_G(q)$, and similarly for other polynomials when the
need arises.

\subsection{Potts correlations and flow counts}
The  Potts correlation functions \eqref{eq:tau}
may be expressed in terms of flow polynomials associated
with a certain Poissonian random graph derived from $G$
by replacing each edge by a random number of copies.

For a vector $\boldm=(m_e: e\in E)$ of non-negative integers,
let $G_\boldm=(V,E_\boldm)$ be the graph with vertex
set $V$ and, for each $e\in E$, with exactly $m_e$ edges
in parallel joining the endvertices of the edge $e$;
the original edge $e$ is  removed.

Let $\bl \ge 0$, and let $\bP=
(\pim_e:e\in E)$ be a family of independent random variables
such that $\pim_e$ has the Poisson 
distribution with parameter $\bl$.
The random graph $G_\bP=(V,E_\bP)$ is
 called a \emph{Poisson graph with intensity $\bl$}. Let
$\Pr_\bl$ and $\EE_\bl$ denote the corresponding probability measure
and expectation.

For $x,y\in V$, 
let $G_\bP^{x,y}$ denote the graph obtained from $G_\bP$ by adding an
edge with endvertices $x$, $y$. If $x$ and $y$ are 
adjacent in the original graph $G_\bP$, we add a
further edge between them. 
Potts correlations are related to flow counts as follows.

\begin{theorem}\label{8.flowcorr} 
Let $q \in\{2,3,\dots\}$ and $\bl \ge 0$. Then
\begin{equation}
q\tau_{\be,q}(x,y) = \frac{\EE_\bl(C(G_\bP^{x,y};q))}
{\EE_\bl(C(G_\bP;q))},
\qq x,y\in V.
\label{8.flowcorr2}
\end{equation}
\end{theorem}

This formula is particularly striking when $q=2$, since
non-zero mod-2 flows take only the value 1. 
A finite graph $H=(W,F)$ is said to be \emph{even} if every vertex 
has even degree.
Evidently $C_H(2)=1$ 
if $H$ is even, and $C_H(2)=0$ otherwise.
By \eqref{8.flowcorr2},
for any graph $G$, 
\begin{equation}
q\tau_{\be,q}(x,y)=\frac{\Pr_\bl(G^{x,y}_\bP\text{ is even})}
{\Pr_\bl(G_\bP\text{ is even})},
\label{8.flowcorr3}
\end{equation}
when $q=2$. This observation is central to the so called
\lq random-current expansion' of the  Ising model, which
has proved very powerful in the study of both classical and quantum Ising models. See \cite{MR723877,MR894398,ADS,MR857063,Bj14}. 

Theorem \ref{8.flowcorr}
may be extended via \eqref{corrconn} to the \rc\ model.
The following is obtained by expressing the flow polynomial
in terms of the Tutte polynomial $T$, and allowing $q$ to 
vary continuously.

\begin{theorem}\label{8.flowconn} 
Let $p\in (0,1)$, $q \in(0,\oo)$, and let $\bl$ satisfy $p=1-e^{-\bl  q}$. 
\begin{romlist}
\item 
For $x,y\in V$,
\begin{equation}
(q-1)\phi_{G,p,q}(x\lra y) = 
\frac {\EE_\bl\bigl((-1)^{1+|E_\bP|} T(G_\bP^{x,y};0,1-q)\bigr)}
{\EE_\bl\bigl((-1)^{|E_\bP|}T(G_\bP;0,1-q)\bigr)}.
\label{8.flowconn2}
\end{equation}
\item
For $q\in\{2,3,\dots\}$,
\begin{equation}
\ZRC(p,q) = \phi_{G,p}(q^{k(\om)})=
\left[(1-p)^{|E|(q-2)/q} q^{|V|}\right] \EE_\bl(C(G_\bP;q)).
\label{8.compflow}
\end{equation}
\end{romlist}
\end{theorem}

Further details may be found in \cite{MR2314565}.

\subsection{The random-current expansion when $q=2$}

Unlike the Potts model, 
there is a fairly complete
analysis of the Ising model. A principal part
in this analysis is played by Theorem \ref{8.flowcorr} with $q=2$
under the heading `random-current expansion'. This has permitted
proofs amongst other things of the exponential decay
of correlations in the low-$\be$ regime on the cubic lattice
$\LL^d$ with $d\ge 2$. It has not so far been
possible to extend this work to general Potts models, but 
Theorem \ref{8.flowcorr}
could play a part in such an extension.

Let $G=(V,E)$ be a finite graph and set $q=2$.
By Theorem \ref{8.flowcorr},
\begin{equation}\label{8.isingeven}
2\tau_{\bl,2}(x,y)=\dfrac{\Pr_\bl(\text{$G_\bP^{x,y}$ is even})}
{\Pr_\bl(\text{$G_\bP$ is even})},
\qq 0\le \bl<\oo.
\end{equation}

There is an important correlation inequality
known as  Simon's inequality, \cite{MR589426}. 
Let $x,z\in V$ be distinct vertices.
A set $W$ of vertices is said to {\it separate\/} $x$ and $z$ if 
$x,z\notin W$ and every path from $x$ to $z$ contains some vertex of $W$.

\begin{theorem}\label{8.simon}
Let $x,z\in V$ be distinct vertices, and let $W$ separate
$x$ and $z$. Then $\kappa_{\bl,2}(x,y):= 2\tau_{\bl,2}(x,y)$ satisfies
$$
\kappa_{\bl,2}(x,z)\le \sum_{y\in W} \kappa_{\bl,2}(x,y)\kappa_{\bl,2}(y,z).
$$
\end{theorem}

The Ising model corresponds to a \rc\ measure $\phi_{p,q}$ with $q=2$.
By \eqref{eq:corrconn},
$$
\kappa_{\bl,q}(x,y)=\phi_{p,q}(x\lra y),
$$
where $p=1-e^{-\bl q}$ and $q=2$. The Simon inequality
may be written in the form
\begin{equation}
\phi_{p,q}(x\lra z)\le\sum_{y\in W}\phi_{p,q}(x\lra y)\phi_{p,q}(y\lra z)
\label{8.simon2}
\end{equation}
whenever $W$ separates $x$ and $z$. It is well known that this inequality
is valid also when $q=1$, see \cite[Chap.\ 6]{MR1707339}. 
One may conjecture that it holds for any $q\in[1,2]$.

\section{The Limit of Zero Temperature}

The physical interpretation of the constant $\be$ is as $\be = 1/(kT)$ where
$k$ is Boltzmann's constant and $T$ denotes (absolute) temperature. The limit
$T\downarrow 0$ corresponds to the limit $\be\to \oo$.  The ferromagnetic
Potts measure $\pi_{\be,q}$ on a finite graph $G=(V,E)$ converges weakly
to the probability measure that allocates a uniform random spin to
each connected component of $G$, this being constant on each component
and independent between components. A realization of this recipe is
called  a `ground state' of the system.

The situation is more interesting in the presence of a vector $\boldj$ of edge-parameters, 
some of which are negative. The ground  states in this case
are colourings $\kappa$ of $V$ with the colour palette $\{1,2,\dots,q\}$ and the property
that, for any each $e = \langle x,y \rangle$,
\begin{equation}\label{eq:frust}
\kappa(x) 
\begin{cases} = \kappa(y) &\text{if } J_e >0,\\
\ne \kappa(y) &\text{if } J_e < 0.
\end{cases}
\end{equation}
In the purely antiferromagnetic case, such a colouring $\kappa$
has the property that any two neighbours have different
colours. There exist graphs $G$ for which \eqref{eq:frust} has no
solution, and such graphs are called `frustrated'.

\begin{theorem}
For the purely antiferromagnetic Potts model,
$$
\ZP(\be,q) \to \chi_G(q) \qq\text{as } \be \to\oo,
$$
where $\chi_G$ is the chromatic polynomial of $G$.
\end{theorem}

\section{The Random-Cluster Model on the Complete Graph}

When the underlying graph is the complete graph $K_n$, the asymptotic
behaviour of the corresponding random-cluster  partition function $\ZRC(n,p,q)$ may be
studied using a mixture of combinatorics and probability, within the regime $q \ge 1$, $p=\la/n$.
Here is some notation and explanation, in preparation for the main theorem.

Let $q \ge 1$ and $p=\la/n$.
It turns out that there is a critical value of $\lambda$ that marks the arrival of a giant
cluster in the \rc\ model on $K_n$, and this value is given by
$$
\lc (q)=\begin{cases}
q&\text{if $q\in(0,2]$,}\\
2\left(\dfrac{q-1}{q-2}\right) \log (q-1)&\text{if $q\in(2,\oo)$}.
\end{cases}
$$
As $\lambda$ increases through the value $\lc$, a giant cluster
of size approximately $\theta(\la,q)n$ is created, where 
$$
\theta (\lambda ,q)=\begin{cases}
0&\text{if $\lambda <\lc (q)$},\\
\tmax&\text{if $\lambda\geq\lc (q)$},
\end{cases}
$$
and $\tmax$ is the largest root of the equation
$$
e^{-\lambda\theta} =\frac{1-\theta}{1+(q-1)\theta}.
$$

\begin{theorem} 
Let $q\in(0,\oo)$ and $\lambda\in(0,\oo)$. We have that 
$$
\frac 1n \log \ZRC(n,\lambda /n,q)\to \eta (\lambda )\qq\text{as
$n\to\infty$},
$$
where 
\begin{align*}
\eta (\lambda )&= \frac{g(\theta (\lambda,q ))}{2q} -\frac{q-1}{2q} \lambda
+\log q,\\
g(\theta )&=-(q-1)(2-\theta )\log (1-\theta )-\bigl[ 2+(q-1)\theta \bigr]\log\bigl[ 1+(q-1)\theta\bigr] .
\end{align*}
\end{theorem}

By Theorem \ref{thm:tutte-rcm}, this provides an asymptotic evaluation of the Tutte polynomial
$T_{n,\lambda/n,q}(u-1,v-1)$ within the quadrant $u,v\in(0,\oo)$.

See \cite{MR1376340} and \cite[Chap.\ 10]{MR2243761} for further details.



\section{Open Problems} 

There is an enormous range of open problems
associated with Ising, Potts, and random-cluster models. 
Of these, there follows a brief selection
some of which may be allied in part to the Tutte polynomial.

\begin{numlist}

\item Prove or disprove some version of negative association
for the uniform forest measure $\USF$ or the uniform connected subgraph
measure $\UCS$. See Conjecture \ref{2.forcs}.

\item Prove or disprove some version of negative association for
the \rc\ measure $\fpq$ with $0<q<1$.

\item Prove or disprove a version  of Simon's inequality for
a \rc\ measure $\fpq$ with $q \in [1,2]$, as in \eqref{8.simon2}.

\item Establish a version of Simon's inequality,
Theorem \ref{8.simon}, for the Potts model with $q \ge 3$.

\item
More generally, find an application  of mod-$q$ flows
to the $q$-state Potts model with $q \ge 3$, as in Theorem \ref{8.flowcorr}.
Develop such an application for real $q>0$, as in Theorem \ref{8.flowconn}.

\item 
A problem of great current interest to probabilists and mathematical
physicists is to understand the geometry of interfaces in the $3$-state Potts
model on the square lattice $\ZZ^2$, and more generally the \rc\ model
on $\ZZ^2$ with $1 <q < 4$. See \cite{Smi10}, and \cite{BDS2}  for a recent reference.

\end{numlist}

\section*{Acknowledgement}
This work was supported in part
by the Engineering and Physical Sciences Research Council under grant EP/103372X/1. 

\bibliographystyle{plain}
\bibliography{tutte-grg-arxiv}

\printindex

\end{document}